\documentclass[12pt]{amsart}

\usepackage{amsmath, amscd, amssymb}
\numberwithin{equation}{section}
\newcommand{\sM}{{\mathcal M}}

\newcommand{\sR}{{\mathcal R}}
\newcommand{\sS}{{\mathcal S}}
\newcommand{\sT}{{\mathcal T}}

\newcommand{\C}{{\mathbb C}}

\newcommand{\N}{{\mathbb N}}

\newcommand{\PP}{{\mathbb P}}
\newcommand{\Q}{{\mathbb Q}}
\newcommand{\R}{{\mathbb R}}

\newcommand{\Z}{{\mathbb Z}}
\newtheorem{theorem}{Theorem}[section]

\newtheorem{lemma}[theorem]{Lemma}

\theoremstyle{remark}

\theoremstyle{definition}

\begin{document}

\title[Localization, Hurwitz Numbers and the Witten Conjecture]
{Localization, Hurwitz Numbers and the Witten Conjecture}

\author{Lin Chen}
\address{Department of Mathematics\\
University of California at Los Angeles}
\email{chenlin@math.ucla.edu}

\author{Yi Li}
\address{Center of Mathematical Sciences\\
Zhejing University\\ Hangzhou, China} \email{yili@cms.zju.edu.cn}

\author{Kefeng Liu}
\address{Center of Mathematical Sciences, Zhejiang University, Hangzhou,
China and Department of Mathematics\\
University of California at Los Angeles} \email{liu@math.ucla.edu,
liu@cms.zju.edu.cn}
\begin{abstract}In this note, we use the method of [3] to
give a simple proof of famous Witten conjecture. Combining the
coefficients derived in our note and this method, we can derive more
recursion formulas of Hodge integrals.
\end{abstract}

\maketitle

\setcounter{tocdepth}{5} \setcounter{page}{1}

\section{Introduction}
The well-known Witten conjecture states that the intersection
theory of the $\psi$ classes on the moduli spaces of Riemann
surfaces is equivalent to the "Hermitian matrix model" of
two-dimensional gravity. All $\psi$-integrals can be efficiently
computed by using the Witten conjecture, first proved by
Kontsevich [6]. For convenience, we use Witten's natation
\begin{equation}
\langle\tau_{\beta_{1}}\cdots\tau_{\beta_{n}}\rangle_{g}:=\int_{\overline{\sM}_{g,n}}\psi^{\beta_{1}}_{1}\cdots\psi^{\beta_{n}}_{n}.
\end{equation}
The natural generating function for the $\psi$-integrals described
above is
\begin{equation}
F_{g}(t):=\sum_{n\geq0}\frac{1}{n!}\sum_{k_{1},\cdots,k_{n}\geq0}t_{k_{1}}\cdots
t_{k_{n}}\langle\tau_{\beta_{1}}\cdots\tau_{\beta_{n}}\rangle_{g},
\ \ \ F(t,\lambda):=\sum_{g\geq0}F_{g}\lambda^{2g-2}.
\end{equation}

For example, the first system of differential equations
conjectured by Witten are the KDV equations [14]. Let
$F(t):=F(t,1)$, define
\begin{equation}
\langle\langle\tau_{\beta_{1}}\cdots\tau_{\beta_{n}}\rangle\rangle:=\frac{\partial}{\partial
t_{k_{1}}}\cdots\frac{\partial}{\partial t_{k_{n}}}F(t),
\end{equation}
then the KDV equations for $F(t)$ are equivalent to a sequence of
recursive relations for $n\geq1$:
\begin{equation}
(2n+1)\langle\langle\tau_{n}\tau^{2}_{0}\rangle\rangle=\langle\langle\tau_{n-1}\tau_{0}\rangle\rangle\langle\langle\tau^{3}_{0}\rangle\rangle+2\langle\langle\tau_{n-1}\tau^{2}_{0}\rangle\rangle\langle\langle\tau^{2}_{0}\rangle\rangle+\frac{1}{4}\langle\langle\tau_{n-1}\tau^{4}_{0}\rangle\rangle.
\end{equation}

In [5] the author gives a simple proof of the Witten conjecture by
first proving a recursion formula conjectured by
 Dijkgraaf-Verlinde-Verlinde in [1], and as corollary they were able to give a simple
proof of the Witten conjecture by using asymptotic analysis. In
this note, we use the method in [3] to prove this recursion
formula in [1], therefore the Witten conjecture without using the
asymptotic analysis. Combining the coefficients derived in our
note and the approach in [3], we can derive more recursion
formulas of Hodge integrals.

\section{Localization and the Hurwitz Numbers}

Let $\mu: \mu_{1}\geq\cdots\geq\mu_{l(\mu)}>0$ be a partition of
$d>0$ and $\overline{\sM}_{g,0}(\PP^{1},\mu)$, the virtual
dimension $r=2g-2+d+l(\mu)$, be the moduli space of relative
stable morphism to $\PP^{1}$. We refer the reader to [10] for the
details in this subsection. Consider the $\C^{*}$-action on
$\PP^{1}$
\begin{equation*}
t\cdot[z^{0}:z^{1}]=[tz^{0}:z^{1}],
\end{equation*}
then it induces an action on $\overline{\sM}_{g,0}(\PP^{1},\mu)$.
There is a branching morphism
\begin{equation}
{\rm Br}: \overline{\sM}_{g,0}(\PP^{1},\mu)\longrightarrow
\PP^{r},
\end{equation}
with this action, the branching morphism is $\C^{*}$-equivariant.
The Hurwitz numbers can be defined by
\begin{equation}
H_{g,\mu}:=\int_{[\overline{\sM}_{g,0}(\PP^{1},\mu)]^{{\rm
virt}}}{\rm Br}^{*}H^{r}
\end{equation}
with $H\in H^{2}(\PP^{r};\Z)$ the hyperplane class.

\subsection{Localization and Hurwitz Numbers}From the localization formula in [8, 10], we
have
\begin{equation}
H_{g,\mu}=(-1)^{k}k!\widetilde{I}^{k}_{g,\mu},
\end{equation}
where $\widetilde{I}^{k}_{g,\mu}$ are the contribution of graphs
of $\overline{\sM}_{g,0}(\PP^{1},\mu)$. The $k=0$ case implies the
well-known ELSV formula [2]:
\begin{equation}
H_{g,\mu}=\frac{r!}{|{\rm
Aut}(\mu)|}\prod^{l(\mu)}_{i=1}\frac{\mu^{\mu_{i}}_{i}}{\mu_{i}!}\int_{\overline{\sM}_{g,l(\mu)}}\frac{\Lambda^{\vee}_{g}(1)}{\prod^{l(\mu)}_{i=1}(1-\mu_{i}\psi_{i})},
\end{equation}
where
$\Lambda^{\vee}_{g}(1)=1-\lambda_{1}+\cdots+(-1)^{g}\lambda_{g}$.
For $k=1$, we derive the cut-and-join equation
\begin{eqnarray}
H_{g,\mu}&=&\sum_{\nu\in J(\mu)}I_{1}(\nu)H_{g,\nu}+\sum_{\nu\in
C(\mu)}I_{2}(\nu)H_{g-1,\nu} \\
&+&\sum_{g_{1}+g_{2}=g}\sum_{\nu^{1}\cup\nu^{2}\in
C(\mu)}\binom{r-1}{2g_{1}-2+|\nu^{1}|+l(\nu^{1})}I_{3}(\nu^{1},\nu^{2})H_{g_{1},\nu^{1}}H_{g_{2},\nu^{2}}.\nonumber
\end{eqnarray}

\subsection{Notations}In this subsection, we explain some
notations appeared in the above subsection. Let $\mu:
\mu_{1}\geq\cdots\geq\mu_{n}>0$, and for each positive integer
$i$,
\begin{equation}
m_{i}(\mu)=|\{j|\mu_{j}=i\}|.
\end{equation}
First, we recall the definitions of $J_{\mu}$ and $C_{\mu}$ (see
[3])
\begin{eqnarray*}
J^{i,j}(\mu)&=&\{(\mu_{1},\cdots,\widehat{\mu}_{i},\cdots,\widehat{\mu}_{j},\cdots,\mu_{n},\mu_{i}+\mu_{j})\},
\ J(\mu) \ = \ \cup^{n}_{i=1}\cup^{n}_{j=i+1}J^{i,j}(\mu); \\
C^{i,p}(\mu)&=&\{(\mu_{1},\cdots,\widehat{\mu}_{i},\cdots,\mu_{n},p,\mu_{i}-p)\},
\ C^{i}(\mu) \ = \ \cup^{\mu_{i}-1}_{p=1}C^{i,p}(\mu), \ C(\mu) \
= \ \cup^{n}_{i=1}C^{i}(\mu).
\end{eqnarray*}
If $\nu\in J^{i,j}(\mu)$, we rewrite $\nu:=\mu^{i,j}$ and then the
$I_{1}(\nu)$ is given by
\begin{equation}
I_{1}(\nu)=\frac{\mu_{i}+\mu_{j}}{1+\delta^{\mu_{i}}_{\mu_{j}}}m_{\mu_{i}+\mu_{j}}(\mu^{i,j})
\end{equation}
Similarly, for $\nu\in C^{i,p}(\mu)$, we can also rewrite
$\nu=\mu^{i,p}$ and the $I_{2}(\nu)$ is defined by
\begin{equation}
I_{2}(\nu)=\frac{p(\mu_{i}-p)}{1+\delta^{p}_{\mu_{i}-p}}m_{p}(\mu^{i,p})(m_{\mu_{i}-p}(\mu^{i,p})-\delta^{p}_{\mu_{i}-p})
\end{equation}
Finally, if $\nu\in C^{i,p}(\mu)$, let $\nu^{1}\cup\nu^{2}=\nu$,
the definition of $I_{3}(\nu^{1},\nu^{2})$ is given by
\begin{equation}
I_{3}(\nu^{1},\nu^{2})=\frac{p(\mu_{i}-p)}{1+\delta^{p}_{\mu_{i}-p}}m_{p}(\nu^{1})m_{\mu_{i}-p}(\nu^{2})
\end{equation}
Define the formal power series
\begin{equation}
\Phi(\lambda,p)=\sum_{\mu}\sum_{g\geq0}H_{g,\mu}\frac{\lambda^{2g-2+|\mu|+l(\mu)}}{(2g-2+|\mu|+l(\mu))!}p_{\mu},
\end{equation}
from [11], we have the following version of cut-and-join equation
\begin{equation}
\frac{\partial\Phi}{\partial\lambda}=\frac{1}{2}\sum_{i.j\geq1}\left(ijp_{i+j}\frac{\partial^{2}\Phi}{\partial
p_{i}\partial p_{j}}+ij p_{i+j}\frac{\partial\Phi}{\partial
p_{i}}\frac{\partial\Phi}{\partial
p_{j}}+(i+j)p_{i}p_{j}\frac{\partial\Phi}{\partial
p_{i+j}}\right).
\end{equation}
At last, we define
\begin{equation}
\Phi_{g,n}(z,p)=\sum_{d\geq1}\sum_{\mu\vdash d,
l(\alpha)=n}\frac{H_{g,\mu}}{r!}p_{\mu}z^{d},
\end{equation}
by simple calculation, we can rewrite (2.12) in the following form
\begin{equation}
\Phi_{g,n}(z;p)=\frac{1}{n!}\sum_{b_{1},\cdots,b_{n}\geq0, 0\leq
k\leq
g}(-1)^{k}\langle\tau_{b_{1}}\cdots\tau_{b_{n}}\lambda_{k}\rangle_{g}\prod^{n}_{i=1}\phi_{b_{i}}(z;p),
\end{equation}
where
\begin{equation}
\phi_{i}(z;p)=\sum_{m\geq0}\frac{m^{m+i}}{m!}p_{m}z^{m}, \ \ \
i\geq0.
\end{equation}

\section{Symmetrization Operator and Rooted Tree Series}

In this section, we use the method in [3] to prove the recursion
formula which implies the Witten conjecture/Kontsevich theorem.
Their method consists of the following steps: (1) introduce three
operators to change the variables; (2) compare the leading
coefficient of both sides to derive the recursion formula which
implies the Witten conjecture. Kim-Liu [5] have proved the Witten
conjecture via the asymptotic analysis which writes each
$\mu_{i}=x_{i}N$ for some $x_{i}\in\R$ and $N\in\N$. The main
problem arising in [2] is the asymptotic expansion of series
\begin{equation*}
e^{-n}\sum_{p+q=n}\frac{p^{p+k+1}q^{q+l+1}}{p!q!}, \ \ \
e^{-n}\sum^{n}_{p+q=n}\frac{p^{p+k+1}q^{q-1}}{p!q!}
\end{equation*}
for any $k, l\in\N$ which are not easy to compute. The idea here
is that by using the method in [3], we can avoid these problems to
derive the recursion formula.

\subsection{Symmetrization Operator $\Box_{n}$}First, we symmetrize $\Phi_{g,n}(z,p)$ by using the linear symmetrization operator $\Box_{n}$
\begin{equation}
\Box_{n}(p_{\alpha}z^{|\alpha|})=\delta_{l(\alpha),n}\sum_{\sigma\in
S_{n}}x^{\alpha_{1}}_{\sigma(1)}\cdots x^{\alpha_{n}}_{\sigma(n)},
\end{equation}
where $S^{n}$ is the $n$-order symmetric group. It is easy to
prove the following lemma via some elementary identities or see
[3].
\begin{lemma}For $n, g\geq1$ or $n\geq3, g\geq0$, then we have
\begin{equation}
\Box_{n}(\Phi_{g,n}(z,p))(x_{1},\cdots,x_{n})=\frac{1}{n!}\sum_{b_{1},\cdots,b_{n}\geq0,
0\leq k\leq
g}(-1)^{k}\langle\tau_{b_{1}}\cdots\tau_{b_{n}}\lambda_{k}\rangle_{g}\sum_{\sigma\in
S_{n}}\prod^{n}_{i=1}\phi_{b_{i}}(x_{\sigma(i)}),
\end{equation}
where
\begin{equation}
\phi_{i}(x):=\phi(x;1)=\sum_{m\geq1}\frac{m^{m+i}}{m!}x^{m}.
\end{equation}
\end{lemma}

\subsection{Rooted Tree Series}The authors introduce the {\it rooted tree
series} $w(x)$ in [3] to simplify the series $\phi(x)$
\begin{equation}
w(x)=\sum_{m\geq1}\frac{m^{m-1}}{m!}x^{m},
\end{equation}
then
\begin{equation}
\phi_{i}(x)=\left(x\frac{d}{dx}\right)^{i+1}w(x):=\nabla^{i+1}_{x}w(x)
\end{equation}
with $\nabla_{x}:=x\frac{d}{dx}$. The rooted tree series is the
unique formal power series solution of the functional equation
(see [3])
\begin{equation}
w(x)=xe^{w(x)}.
\end{equation}

Let $y(x):=\frac{1}{1-w(x)}$ and $y_{j}=y(x_{j})$, we consider
change of variables using the operator
\begin{equation}
L: \Q[[x_{1},\cdots,
x_{n}]]\longrightarrow\Q[[y_{1},\cdots,y_{n}]], \ \ \
f(x_{1},\cdots,x_{n})\longmapsto f(y_{1},\cdots,y_{n}).
\end{equation}

\begin{lemma}Denote $w_{j}=w(x_{j})$, then
\begin{eqnarray}
L\nabla_{x_{j}}&=&(y^{2}_{j}-y_{j})\nabla_{y_{j}}L, \ \ \
L\nabla_{w_{j}} \ = \ (y_{j}-1)\nabla_{y_{j}}L, \\
L(\phi_{i}(x_{j}))&=&[(y^{2}_{j}-y_{j})\nabla_{y_{j}}]^{i}(y_{j}-1),
\ \ i\geq0.
\end{eqnarray}
\end{lemma}
\begin{proof}Differentiating the functional equation (2.6), we obtain
\begin{equation*}
\nabla_{x_{j}}=\frac{1}{1-w_{j}}\nabla_{w_{j}}.
\end{equation*}
Note that $dy_{j}=y^{2}_{j}dw$, then
$\nabla_{w_{j}}=(y_{j}-1)\nabla_{y_{j}}$ and
\begin{equation*}
L\nabla_{x_{j}}=L\left(\frac{1}{1-w_{j}}\nabla_{w_{j}}\right)=y_{j}(y_{j}-1)\nabla_{y_{j}}L.
\end{equation*}
For the similar reason, the reader can prove the rest identities.
\end{proof}

\subsection{The Coefficients of $L(\phi_{i}(x_{j}))$}Let
\begin{equation}
F_{i}(y):=[(y^{2}-y)\nabla_{y}]^{i}(y-1):=\sum^{i+1}_{j=1}f(j,i)y^{2i+2-j},
\end{equation}
then from (3.9) we get
\begin{eqnarray*}
\sum^{i+2}_{j=1}f(j,i+1)y^{2i+4-j}&=&(2i+1)f(1,i)y^{2i+3}-if(i+1,i)y^{i+2}
\\
&+&\sum^{i+1}_{j=2}[(2i+2-j)f(j,i)-(2i+3-j)f(j-1,i)]y^{2i+4-j}.
\end{eqnarray*}
Comparing the coefficients of both sides equation, we have the
crucial identities
\begin{eqnarray*}
f(1,i+1)&=&(2i+1)f(1,i), \\
f(j,i+1)&=&(2i+2-j)f(j,i)-(2i+3-j)f(j-1,i), \ \ \ 2\leq j\leq i+2 \\
f(i+1,i+1)&=&-if(i+1,i).
\end{eqnarray*}
The initial value of the coefficients are given by
\begin{equation}
f(1,1)=1, \ \ \ f(2,1)=-1,
\end{equation}
hence
\begin{equation}
f(1,i)=(2i-1)!!, \ \ \ f(i+1,i)=(-1)^{i}i!.
\end{equation}
Now, we only give the explicit expression of $f(2,i)$ while the
similar approach can be used to derive other explicit expressions.
If we write
\begin{equation}
a(j,i)=\frac{f(j,i)}{(2i-j)!!}, \ \ \ 2\leq j\leq i+1,
\end{equation}
note that $a(j,1)=-1$, we have a linear equation system
\begin{equation*}
a(j,i+1)=a(j,i)-\frac{2i+3-j}{(2i+2-j)!!}f(j-1,i).
\end{equation*}
Finally, the recursion formula of $f(j,i)$ is given by
\begin{equation}
f(j,i)=-(2i-j)!!\left[1+\sum^{i-1}_{k=1}\frac{2k+3-j}{(2k+2-j)!!}f(j-1,k)\right],
\ \ 2\leq j\leq i+1, i\geq1.
\end{equation}
For $i=2$, it turns to the explicit expression
\begin{equation}
f(2,i)=-(2i-2)!!\left[1+\sum^{i-1}_{k=1}\frac{(2k+1)!!}{(2k)!!}\right].
\end{equation}

\begin{lemma}We can rewrite $f(2,i)$ as
\begin{equation}
f(2,i)=-\frac{(2i+1)!!}{3},
\end{equation}
and
\begin{equation*}
2i\cdot f(2,i)=-\frac{(2i+1)!}{2^{i-1}(i-1)!}.
\end{equation*}
\end{lemma}
\begin{proof}By induction on $i$, we obtain
\begin{equation*}
f(2,i+1)=-\frac{(2i)(2i+1)!!+3(2i+3)!!}{3}=-\frac{(2i+3)!!}{3}.
\end{equation*}
\end{proof}

\section{Proof of the Dijkgraaf-Verlinde-Verlinde Conjecture}

For $i,j\geq 0$, $i+j\leq n$, let
$\underset{i,j}{\overset{x}{\Box}}$ be the mapping, applied to a
series in $x_1, \cdots ,x_n$, given by
\begin{equation}
\underset{i,j}{\overset{x}{\Box}}f(x_1, \cdots
,x_n)=\sum_{\sR,\sS,\sT}f(x_{\sR},x_{\sS},x_{\sT}),
\end{equation}
where the summation is over all ordered partitions $(\sR,\sS,\sT)$
of $\{1, \cdots, n\}$, where $\sR=\{x_{r_1}, \cdots, x_{r_i}\}$,
$\sS=\{x_{s_1}, \cdots, x_{s_j}\},$ $\sT=\{x_{t_1}, \cdots,
x_{t_{n-i-j}}\}$ and $$(x_{\sR},x_{\sS},x_{\sT})=(x_{r_1}, \cdots,
x_{r_i},x_{s_1}, \cdots ,x_{s_j},x_{t_1},\cdots ,x_{t_{n-i-j}}),$$
and where $r_1< \cdots <r_i$, $s_1< \cdots <s_j$, and $t_1< \cdots
<t_{n-i-j}.$ The following result gives an expression for the
result of applying the symmetrization operator $\Box_{n}$ to the
cut-and-join equation for $\Phi_{g,n}(z,p)$. Denote
$\triangle_{y_{j}}:=(y^{2}_{j}-y_{j})\nabla_{y_{j}}$. Applying the
symmetrization operator $\Box_{n}$ to the join-cut Equation, [3]
have proved the following version of join-and-cut equation
\begin{equation}
\left(\sum^{n}_{i=1}(y_{i}-1)\nabla_{y_{i}}+n+2g-2\right)L\Box_{n}\Phi_{g,n}(y_{1},\cdots,y_{n})=T'_{1}+T'_{2}+T'_{3}+T'_{4},
\end{equation}
where
\begin{eqnarray*}
T'_{1}&=&\frac{1}{2}\sum^{n}_{i=1}\left(\triangle_{y_{i}}\triangle_{y_{n+1}}L\Box_{n+1}\Phi_{g-1,n+1}(y_{1},\cdots,y_{n+1})\right)|_{y_{n+1}=y_{i}},
\\
T'_{2}&=&\underset{1,1}{\overset{y}{\Box}}y^{2}_{1}\frac{y_{2}-1}{y_{1}-y_{2}}\triangle_{y_{1}}L\Box_{n-1}\Phi_{g,n-1}(y_{1},y_{3},\cdots,y_{n}),
\\
T'_{3}&=&\sum^{n}_{k=3}\underset{1,k-1}{\overset{y}{\Box}}\left(\triangle_{y_{1}}L\Box_{k}\Phi_{0,k}(y_{1},\cdots,y_{k})\right)\left(\triangle_{y_{1}}L\Box_{n-k+1}\Phi_{g,n-k+1}(y_{1},y_{k+1},\cdots,y_{n})\right),
\\
T'_{4}&=&\frac{1}{2}\sum_{1\leq k\leq n, 1\leq a\leq
g-1}\underset{1,k-1}{\overset{y}{\Box}}\left(\triangle_{y_{1}}L\Box_{k}\Phi_{a,k}(y_{1},\cdots,y_{k})\right)\left(\triangle_{y_{1}}L\Box_{n-k+1}\Phi_{g-a,n-k+1}(y_{1},y_{k+1},\cdots,y_{n})\right).
\end{eqnarray*}

\subsection{Expansions}First we have the
following expansion formula
\begin{equation}
L\left(\prod^{n}_{i=1}\phi_{b_{i}}(x_{\sigma(i)})\right)=\prod^{n}_{i=1}(2b_{i}-1)!!y^{2b_{i}+1}_{\sigma(i)}+{\rm
lower \ terms}.
\end{equation}
From this point, we see that the polynomial
$L\Box_{n}H^{g}_{n}(y_{1},\cdots,y_{n})$ can be written
as
\begin{equation*}
L\Box_{n}\Phi_{g,n}(y_{1},\cdots,y_{n})=\sum_{b_{1}+\cdots+b_{n}=3g-3+n}\langle\tau_{b_{1}}\cdots\tau_{b_{n}}\rangle_{g}\prod^{n}_{i=1}(2b_{i}-1)!!y^{2b_{i}+1}_{i}+{\rm
l.t.},
\end{equation*}
where l.t. means lower order terms. We write the left hand side of
equation (4.2) by {\it LHS} while another side by ${RHS}_{1},
{RHS}_{2}, {RHS}_{3}$ and ${RHS}_{4}$, then
\begin{eqnarray*}
LHS&=&\sum^{n}_{i=1}y_{i}\nabla_{y_{i}}\sum_{b_{1}+\cdots+b_{n}=3g-3+n}\left[\langle\tau_{b_{1}}\cdots\tau_{b_{n}}\rangle_{g}(2b_{1}-1)!!\cdots(2b_{n}-1)!!\right]\prod^{n}_{l=1}y^{2b_{l}+1}_{l}+{\rm
l.t.} \\
&=&\sum_{b_{1}+\cdots+b_{n}=3g-3+n}\left[\langle\tau_{b_{1}}\cdots\tau_{b_{n}}\rangle_{g}(2b_{1}-1)!!\cdots(2b_{n}-1)!!\right]\sum^{n}_{i=1}(2b_{i}+1)y_{i}\prod^{n}_{l=1}y^{2b_{l}+1}_{l}+{\rm
l.t.} \\
RHS_{1}&=&\frac{1}{2}\sum_{b_{1}+\cdots+b_{n+1}=3g-5+n}\left[\langle\tau_{b_{1}}\cdots\tau_{b_{n+1}}\rangle_{g-1}(2b_{1}-1)!!\cdots(2b_{n+1}-1)!!\right]
\\
&\cdot&\sum^{n}_{i=1}\left((2b_{i}+1)(2b_{n+1}+1)y^{2}_{i}y^{2}_{n+1}\prod^{n+1}_{l=1}y^{2b_{l}+1}_{l}\right)\left.\right|_{y_{i}=y_{n+1}}+{\rm
l.t.} \\
RHS_{2}&=&\underset{1,1}{\overset{y}{\Box}}\left(\sum_{b_{1}+b_{3}+\cdots+b_{n}=3g-4+n}\left[(2b_{1}+1)!!(2b_{3}-1)!!\cdots(2b_{n}-1)!!\langle\tau_{b_{1}}\tau_{b_{3}}\cdots\tau_{b_{n}}\rangle_{g}\right]\right.
\\
&\cdot&\left.\sum_{m\geq0}\left(\frac{y_{2}}{y_{1}}\right)^{m}y_{2}y^{3}_{1}y^{2b_{1}+1}_{1}\prod^{n}_{l=2}y^{2b_{l}+1}_{l}\right)
+{\rm l.t.} \\
RHS_{3}&=&\sum^{n}_{k=3}\underset{1,k-1}{\overset{y}{\Box}}\left(\sum_{b_{1}+\cdots+b_{k}=k-3}(2b_{1}-1)!!\cdots(2b_{k}-1)!!\langle\tau_{b_{1}}\cdots\tau_{b_{k}}\rangle_{0}(2b_{1}+1)y^{2}_{1}\prod^{k}_{l=1}y^{2b_{l}+1}_{l}\right) \\
&\cdot&\left(\sum_{\overline{b}_{1}+b_{k+1}+\cdots+b_{n}=3g-k-2+n}\left[(2\overline{b}_{1}-1)!!(2b_{k+1}-1)!!\cdots(2b_{n}-1)!!\langle\tau_{\overline{b}_{1}}\tau_{b_{k+1}}\cdots\tau_{b_{n}}\rangle_{g}\right]\right.
\\
&\cdot&\left.(2\overline{b}_{1}+1)y^{2}_{1}\prod^{n}_{l=k+1}y^{2b_{l}+1}_{l}y^{2\overline{b}_{1}+1}_{1}\right)+{\rm
l.t.} \\
RHS_{4}&=&\frac{1}{2}\sum_{1\leq k\leq n, 1\leq a\leq
g-1}\underset{1,k-1}{\overset{y}{\Box}} \\
&\cdot&\left(\sum_{b_{1}+\cdots+b_{k}=3a-3+k}(2b_{1}-1)!!\cdots(2b_{k}-1)!!\langle\tau_{b_{1}}\cdots\tau_{b_{k}}\rangle_{a}(2b_{1}+1)y^{2}_{1}\prod^{k}_{l=1}y^{2b_{l}+1}_{l}\right) \\
&\cdot&\left(\sum_{\overline{b}_{1}+b_{k+1}+\cdots+b_{n}=3g-k-2+n-3a}\left[(2\overline{b}_{1}-1)!!(2b_{k+1}-1)!!\cdots(2b_{n}-1)!!\langle\tau_{\overline{b}_{1}}\tau_{b_{k+1}}\cdots\tau_{b_{n}}\rangle_{g-a}\right]\right.
\\
&\cdot&\left.(2\overline{b}_{1}+1)y^{2}_{1}\prod^{n}_{l=k+1}y^{2b_{l}+1}_{l}y^{2\overline{b}_{1}+1}_{1}\right)+{\rm
l.t.}
\end{eqnarray*}

\subsection{Picking Out the Coefficients}Now, we can only consider the
coefficients of monomial $y^{2(b_{1}+1)}_{1}y^{2b_{2}+1}_{2}\cdots
y^{2b_{n}+1}_{n}$ ($b_{1}+\cdots+b_{n}=3g-3+n$) on both sides of
equation (4.2). By simply calculating, we find
\begin{eqnarray*}
LHS&=&(2b_{1}+1)!!(2b_{2}-1)!!\cdots(2b_{n}-1)!!\langle\tau_{b_{1}}\cdots\tau_{b_{n}}\rangle_{g}
\\
RHS_{1}&=&\frac{1}{2}\sum_{a+b=b_{1}-2}(2a+1)!!(2b+1)!!\prod^{n}_{l=2}(2b_{l}-1)!!\langle\tau_{a}\tau_{b}\tau_{b_{2}}\cdots\tau_{b_{n}}\rangle_{g-1}
\\
RHS_{2}&=&\sum^{n}_{l=2}(2(b_{1}+b_{l}-1)+1)!!(2b_{2}-1)!!\cdots(2b_{l-1}-1)!!(2b_{l+1}-1)!!\cdots(2b_{n}-1)!!
\\
&\cdot&\langle\sigma_{b_{1}+b_{l}-1}\sigma_{b_{2}}\cdots\sigma_{b_{l-1}}\sigma_{b_{l+1}}\cdots\sigma_{b_{n}}\rangle_{g}
\\
RTS_{3,4}&=&\frac{1}{2}\sum_{X\cup
Y=S}\sum_{a+b=b_{1}-2}\sum_{g_{1}+g_{2}=g}(2a+1)!!(2b+1)!!\prod^{n}_{l=2}(2b_{l}-1)!!\langle\tau_{a}\prod_{\alpha\in
X}\tau_{\alpha}\rangle_{g_{1}}\langle\tau_{b}\prod_{\beta\in
Y}\tau_{\beta}\rangle_{g_{2}},
\end{eqnarray*}
where $S=\{b_{2},\cdots, b_{n}\}$. Finally, multiplying the constant
$(2b_{2}+1)\cdots(2b_{n}+1)$, we obtain the recursion formula as
conjectured by Dijkgraaf-Verlinde-Verlinde which implies the Witten
conjecture
\begin{eqnarray*}
\langle\widetilde{\tau}_{b_{1}}\prod^{n}_{l=2}\widetilde{\tau}_{b_{l}}\rangle_{g}&=&\sum^{n}_{l=2}(2b_{l}+1)\langle\widetilde{\tau}_{b_{1}+b_{l}-1}\prod^{n}_{k=2,
k\neq
l}\widetilde{\tau}_{b_{k}}\rangle_{g}+\frac{1}{2}\sum_{a+b=b_{1}-2}\langle\widetilde{\tau}_{a}\widetilde{\tau}_{b}\prod^{n}_{l=2}\widetilde{\tau}_{b_{l}}\rangle_{g-1}
\\
& &\frac{1}{2}\sum_{X\cup Y=\{b_{2},\cdots,
b_{n}\}}\sum_{a+b=b_{1}-2,
g_{1}+g_{2}=g}\langle\widetilde{\tau}_{a}\prod_{\alpha\in
X}\widetilde{\tau}_{\alpha}\rangle_{g_{1}}\langle\widetilde{\tau}_{b}\prod_{\beta\in
Y}\widetilde{\tau}_{\beta}\rangle_{g_{2}}.
\end{eqnarray*}
where $\widetilde{\tau}_{b_{l}}=[(2b_{l}+1)!!]\tau_{b_{l}}$.

\end{document}